# Sistemas $C-$ortocéntricos, bisectrices y euclidianidad en planos de Minkowski

($C-$orthocentric systems, angular bisectors and euclidianity in Minkowski planes)

Tobías de Jesús Rosas Soto
(tjrosas@gmail.com; trosas@demat-fecluz.org)

Departamento de Matemática
Facultad Experimental de Ciencias
Universidad del Zulia
Venezuela

**Resumen**

Mediante el estudio de propiedades geométricas de los sistemas $C-$ortocéntricos, relacionadas con las nociones de ortogonalidad (Birkhoff, isósceles, cordal), bisectriz (Busemann, Glogovskij) y línea soporte a una circunferencia, se muestran nueve caracterizaciones de euclidianidad para planos de Minkowski arbitrarios. Tres de estas generalizan caracterizaciones dadas para planos de Minkowski estrictamente convexos en [8, 9], y las otras seis son nuevos aportes sobre el tema.

**Palabras claves:** Sistemas $C-$ortocéntricos, planos de Minkowski, ortogonalidad, bisectriz, líneas soporte, euclidianidad.

**Abstract**

By studying geometric properties of $C-$orthocentric systems related to the notions of orthogonality (Birkhoff, isosceles, chordal), angular bisector (Busemann, Glogovskij) and support line to a circumference, nine characterizations of the Euclidean plane are shown for arbitrary Minkowski planes. Three of these generalized characterizations given for strictly convex Minkowski planes in [8, 9], and the other six are new contributions on the subject.

## 1. Introducción

En el 2010, los matemáticos S. Wu y H. Martini lograron establecer una serie de caracterizaciones de euclidianidad para planos de Minkowski estrictamente convexos, mediante el estudio de ciertas propiedades geométricas de los sistemas $C-$ortocéntricos en dichos planos (ver [9]). Ellos demostraron que al tomar dos elementos en el plano, ortogonales isósceles, podían construir un sistema $C$–ortocéntrico que cumplía ciertas propiedades geométricas relacionadas con las nociones de ortogonalidad Birkhoff e isósceles, bisectriz de Busemann, líneas soporte, triángulos isósceles y medianas. En estos resultados se basaron para mostrar ciertas caracterizaciones de euclidianidad. Además, ellos usaron fuertemente la condición de convexidad estricta del plano, pues ésta garantizaba algunas condiciones tales como la unicidad en las nociones de ortogonalidad isósceles y Birkhoff.

En el presente año (2014), los matemáticos T. Rosas y W. Pacheco mostraron en [15] una nueva forma de definir la noción de $C$–ortocentro para planos de Minkowski en general, logrando prescindir de la condición de convexidad estricta del plano. Esto les permitió demostrar que al tener dos elementos en un plano de Minkowski arbitrario, ortogonales isósceles, también se podía construir un sistema $C$–ortocéntrico que cumpliera ciertas propiedades geométricas relacionadas con las nociones de ortogonalidad Birkhoff e isósceles, bisectriz de Busemann, triángulos isósceles y medianas.

Continuando con estas ideas, en el presente artículo se estudian algunas propiedades geométricas de los sistemas $C$–ortocéntricos relacionadas con las nociones de ortogonalidad Birkhoff e isósceles, bisectriz de Busemann y Glogovskij, líneas soporte, triángulos isósceles y medianas. A través de estas se presentan ocho caracterizaciones de euclidianidad para planos de Minkowski en general: tres generalizan resultados equivalentes en planos de Minkowski estrictamente convexos (ver [8, 9]), y las seis restantes son nuevos aportes sobre el tema de caracterizaciones. Para estudiar la geometría en planos de Minkowski y las propiedades básicas de las ortogonalidades isósceles, Birkhoff y cordal véanse las referencias [6, 7, 8, 10, 11, 14, 15] y la monografía [1].

Denotemos por $\left(\mathbf{R}^2, \|\circ\|\right) = M$ a un plano de Minkowski cualquiera con origen $O$, circunferencia unitaria $C$ y norma $\|\circ\|$, y $M^*$ su espacio dual. Para cualquier punto $x \in M$ y $\lambda \in \mathbf{R}^+$, llamemos al conjunto $C(x, \lambda) = x + \lambda C$ la circunferencia de centro en $x$ y radio $\lambda$. Para $x \neq y$, denotemos por $\langle x, y \rangle$ a la línea que pasa por $x$ y $y$, por $[x, y]$ el segmento entre $x$ y $y$, y por $[x, y\rangle$ el rayo que inicia en $x$ y pasa a través de $y$ (ver [9, 14, 15, 16]).

Dado un punto $x$ y la $C(x,\beta)$ en $M$, si tomamos dos puntos diferentes $v,w \in C(x,\beta)$, entonces la línea $\langle v,w \rangle$ divide el plano en dos semiplanos, $L_+$ y $L_-$. Si el segmento $[y,w]$ no está contenido en $C(x,\beta)$, entonces la línea $\langle v,w \rangle$ divide a $C(x,\beta)$ en los arcos $C(x,\beta) \cap L_+$ y $C(x,\beta) \cap L_-$, entre los puntos $v$ y $w$, los cuales denotaremos por $Arc^+_{C(x,\beta)}(v,w)$ y $Arc^-_{C(x,\beta)}(v,w)$, respectivamente. Luego, si $\langle v,w \rangle$ está contenido en $C(x,\beta)$, entonces $Arc^+_{C(x,\beta)}(v,w)$ es exactamente el segmento $[y,w]$ y $Arc^-_{C(x,\beta)}(v,w)$ es el conjunto $\{v,w\} \cup \{C(x,\beta): s \notin [v,w]\}$ (ver [16]).

Para puntos $x,y \in M$, se dice que $x$ es ortogonal isósceles a $y$ si $\|x+y\| = \|x-y\|$, y lo denotaremos por $x \perp_I y$. De igual forma, $x$ se dice ortogonal Birkhoff a $y$ si $\|x+ty\| \geq \|x\|$ para todo $t \in \mathbf{R}$, y lo denotaremos por $x \perp_B y$ (ver [10, 11]). Por otra parte, sean $L_1 = [p_1,q_1]$ y $L_2 = [p_2,q_2]$ cuerdas de una circunferencia $G$ en $M$, entonces $L_1$ y $L_2$ se dicen ortogonales cordales si la línea que pasa a través de $q_2$ y $p'_2$ (el puesto de $p_2$ en $G$) es paralela a la línea $\langle p_1,q_1 \rangle$. En caso de que $p'_2 = q_2$, se dirá que $L_1$ es ortogonal cordal a $L_2$ si existe una línea soporte de $G$ que pase por $q_2$ y sea paralela a la línea $\langle p_1,q_1 \rangle$. Si $p'_2$ es el opuesto de $q_2$ en $G$, entonces $p'_2 q_2 p_2 q'_2$ es un paralelogramo (ver [8]). Consideraremos la ortogonalidad cordal de forma natural solo con respecto a la circunferencia unitaria $C$, usando el símbolo $\perp_C$ para esto, es decir, escribir $[p_1,q_1] \perp_C [p_2,q_2]$ automáticamente presupone que $[p_1,q_1]$ y $[p_2,q_2]$ son cuerdas de $C$.

Algunas propiedades de estas ortogonalidades que serán de utilidad son las siguientes (ver [8, 10, 11]):

- Para todo $x \in C$ existe un $y \in C$ tal que $x \perp_I y$.
- Para $x,y \in M$ tales que $x \perp_I y$, entonces $y \perp_I x$, es decir, la ortogonalidad isósceles es simétrica.
- Para todo $x \in C$ existe un $y \in C$ tal que $x \perp_B y$.
- Para $x,y \in M$ tales que $x \perp_B y$, entonces $kx \perp_B my$ para todo $k,m \in \mathbf{R}$, es decir, la ortogonalidad Birkhoff es homogénea.
- Para toda cuerda $[x,y] \in C$ existe una cuerda $[z,w] \in C$ tal que $[x,y] \perp_C [z,w]$.

Para rayos no colineales $[p,a\rangle$ y $[p,b\rangle$, llamaremos *ángulo* $\angle apb$ a la capsula convexa formada por los rayos $[p,a\rangle$ y $[p,b\rangle$ (ver [14]). El rayo

$$\left[p, \frac{1}{2}\left(\frac{a-p}{\|a-p\|} + \frac{b-p}{\|b-p\|}\right) + p\right\rangle$$

es llamado la *bisectriz de Busemann* del ángulo $\angle apb$ y se denota por $A_B([p,a\rangle, [p,b\rangle)$ (ver [3]). Notemos que cuando $\|a-p\| = \|b-p\|$, entonces

$$A_B([p,a\rangle, [p,b\rangle) = \left[p, \frac{a+b}{2}\right\rangle.$$

La *bisectriz de Glogovskij* del ángulo $\angle apb$, se define como el conjunto de puntos que equidistan de los rayos $[p,a\rangle$ y $[p,b\rangle$, y se denota por $A_G([p,a\rangle, [p,b\rangle)$ (ver [17]).

Dado que los espacios $M$ y $M^*$ los podemos identificar vía isomorfismos con $\mathbf{R}^2$, entonces podemos ver a $M^*$ como $(M, \|\circ\|^*)$, con $\|\circ\|^*$ la norma dual. Como para $x \in M$ y una recta $L$ en $M$ tal que $x \notin L$, existe un funcional $f$ en $M^*$ tal que $f$ anula a $L$ y $d(x, L) = \frac{f(x)}{\|f\|^*}$ (Teorema de Hanh-Banach), entonces la bisectriz de Glogovskij del ángulo $\angle apb$ se puede expresar como $A_G([p,a\rangle, [p,b\rangle) = \{x \in M : d(x, [p,a\rangle) = d(x, [p,b\rangle)\}$, con

$$d(x, [p,a\rangle) = d(x, f) = \frac{|f(x)|}{\|f\|^*},$$

donde $f$ es el *funcional* asociado a la recta $\langle p, a\rangle$. Además, denotaremos por $A_G^*([p,a\rangle, [p,b\rangle)$ a la bisectriz de Glogovskij del $\angle apb$ en $M^* = (M, \|\circ\|^*)$ (ver [16]).

Si $M$ fuera un plano euclídeo con $w \in M$ y $C(x, \lambda)$ una circunferencia en $M$, la *potencia de $w$ con respecto a la circunferencia* $C(x, \lambda)$ es $\|x - w\|^2 - \lambda^2$. Si $C(y, \alpha)$ es otra circunferencia en $M$ no concéntrica con $C(x, \lambda)$, se llama *eje radical* al lugar geométrico de los puntos que tienen igual potencia con respecto de las circunferencias $C(x, \lambda)$ y $C(y, \alpha)$, es decir, $\{a \in M : \|x - a\|^2 - \lambda^2 = \|x - a\|^2 - \alpha^2\}$ (ver [12]).

Para $p \in M$, denotemos por $S_p$ a la simetría con respecto al punto $p$, dada por la expresión $S_p(w) = 2p - w$ para $w \in M$. $H_{p,k}$ denotará la homotecia con centro $p$ y razón $k$ ($k \in \mathbf{R}$), y está dada por $H_{p,k}(w) = (1-k)p + kw$ para $w \in M$. Recordemos que las simetrías son isometrías en planos de Minkowski, es decir, $\|S_p(w) - S_p(v)\| = \|w - v\|$ para todo $w, v \in M$ (ver [15, 16]).

Diremos que una recta $L$ en $M$, con vector director $v$, es la línea soporte de la circunferencia $C(x,\lambda)$, si para todo $e \in L$ y todo $t \in \mathbf{R}$ se cumple que $\|(x-e)+tv\| \geq \lambda$, es decir, una recta $L$ se dirá línea soporte de la circunferencia $C(x,\lambda)$ en el punto $e$, si $e \in L \cap C(x,\lambda)$ y $d(x,L) = \lambda$ (ver [14, 15, 16]).

Para puntos $x_1, x_2, x_3 \in M$ denotemos por $\Delta x_1 x_2 x_3$ al triángulo de vértices $x_1, x_2, x_3$. Diremos que $p$ es un circuncentro del $\Delta x_1 x_2 x_3$, si $\|p - x_1\| = \|p - x_2\| = \|p - x_3\|$. Denotemos por $C(\Delta x_1 x_2 x_3)$ el conjunto de los circuncentros del $\Delta x_1 x_2 x_3$. Dado un punto $p_4$ y un $\Delta x_1 x_2 x_3$ en $M$, diremos que el $\Delta p_1 p_2 p_3$ es el $p_4$-antitriángulo del $\Delta x_1 x_2 x_3$, si $p_i = S_{m_i}(p_4)$ para $i = 1,2,3$, con $m_i$ los puntos medios de los lados del $\Delta x_1 x_2 x_3$. Si $p_4 \in C(\Delta x_1 x_2 x_3)$, diremos que $x_4$ es el $C$–ortocentro del $\Delta x_1 x_2 x_3$ asociado a $p_4$ si $S_q(p_4) = x_4$, donde $q$ es el punto de simetría del $\Delta x_1 x_2 x_3$ y su $p_4$-antitriángulo. Además, el conjunto $\{x_1, x_2, x_3, x_4\}$ se dice un sistema $C$–ortocéntrico (ver [15, 16]).

## 2. Preliminares.

Los siguientes resultados son necesarios para la investigación.

**Lema 2.1:** Las siguientes proposiciones son válidas:
1. Si para cualquier $x, y \in M$, con $x \perp_I y$, existe un numero $t \in \mathbf{R}^+ - \{1\}$ tal que $x \perp_I ty$, entonces $M$ es euclidiano (ver [2, 9]).
2. Un plano de Minkowski $M$ es euclidiano si y solo si para todo $x, y \in M$ se cumple que si $x \perp_I y$, entonces $x \perp_B y$ (ver [2]).
3. Un plano de Minkowski $M$ es euclidiano si y solo si para todo $x, y \in M$ se cumple que si $x \perp_B y$, entonces $x \perp_I y$ (ver [2]).
4. Un plano de Minkowski es euclidiano si y solo si los puntos medios de cualquier familia de cuerdas paralelas están en una misma recta (ver [2]).

**Teorema 2.1:** (ver [15, 16]) Sea $M$ un plano de Minkowski. Sean $x_1, x_2, x_3$ y $p_4$ puntos en $M$. Sean $m_1, m_2$ y $m_3$ los puntos medios de los segmentos $[x_2, x_3]$, $[x_1, x_3]$ y $[x_1, x_2]$, respectivamente. Definamos los puntos $p_i = S_{m_i}(p_4)$, para $i = 1,2,3$, entonces se cumple lo siguiente:
1. Los segmentos $[x_i, p_i]$ tienen el mismo punto medio $q$, para $i = 1,2,3$. Además,
   $2(q - m_i) = x_i - p_4$ para $i = 1,2,3$, es decir, $q = \dfrac{x_1 + x_2 + x_3 - p_4}{2}$.
2. Si $x_4 = S_q(p_4)$, entonces $x_i - x_j = p_j - p_i$ para $\{i, j\} \subset \{1,2,3,4\}$.
3. $x_i - p_j = p_k - x_l$, donde $\{i, j, k,\} = \{1,2,3,4\}$.

4. Si $g = \dfrac{x_1 + x_2 + x_3}{3}$, entonces $H_{g,-2}(p_4) = x_4$.

**Lema 2.2:** (ver [15, 16]) Sea $M$ un plano de Minkowski con origen $O$. Para cualesquiera $x, z \in M$ con $x \perp_I z$, sea $p_3 = -z$, $p_4 = z$, $x_1 = x$, $x_2 = -x$ y $\lambda = \|x + z\|$. Entonces existen puntos $x_3 \in C(p_4, \lambda)$ y $q \in C(O, \tfrac{\lambda}{2})$ tales que, los conjuntos $\{p_1, p_2, p_3, p_4\}$ y $\{x_1, x_2, x_3, x_4\}$ son sistemas $C$–ortocéntricos, donde $p_1 = S_q(x_1)$, $p_2 = S_q(x_2)$ y $x_4 = S_q(p_4)$. Además, se cumple uno de los siguientes enunciados:

1. $p_1 \in C(x_2, \lambda)$, $p_2 \in C(x_1, \lambda)$ y $x_4 \in C(p_3, \lambda)$.
2. Existen $p_1$ y $p_2$ tales que $\|p_3 - p_1\| = \|p_3 - p_2\|$.
3. $p_1 = \pm\dfrac{\lambda z}{\|z\|} - x$ y $p_2 = \pm\dfrac{\lambda z}{\|z\|} + x$ si y solo si $p_4 \in \left\langle p_3, \dfrac{p_1 + p_2}{2} \right\rangle$. En particular, $p_1 = \dfrac{\lambda z}{\|z\|} - x$ y $p_2 = \dfrac{\lambda z}{\|z\|} + x$ si y solo si $p_4 \in \left[ p_3, \dfrac{p_1 + p_2}{2} \right]$.
4. Si $\|z\| < \lambda$ y $q \in Arc^+_{C(o,\tfrac{\lambda}{2})}\left(\dfrac{p_4 + x_2}{2}, \dfrac{p_4 + x_1}{2}\right)$, con $p_1 \neq p_4, S_{x_2}(p_3)$, entonces $p_3$ y la línea $\langle p_1, p_2 \rangle$ están separados por $L_1$.
5. Existen puntos $p_1$ y $p_2$ tal que $p_3$ y la línea $\langle p_1, p_2 \rangle$ están separados por $L_1$, o $L_1 = \langle p_1, p_2 \rangle$. Además, $p_4 \in A_B([p_3, p_1\rangle, [p_3, p_2\rangle)$.

donde $L_1 = \langle S_{x_1}(p_3), S_{x_2}(p_3) \rangle$.

## 3. Resultados y pruebas.

El siguiente lema nos da las herramientas para probar la mayoría de los resultados sobre caracterización de euclidianidad mencionadas, mediante el estudio de propiedades geométricas en los sistemas $C$ – ortocéntricos en planos normados.

**Lema 3.1:** Sea $M$ un plano de Minkowski, con origen $O$. Para cualesquiera $x, z \in M$ con $x \perp_I z$. Sean $p_3 = -z$, $p_4 = z$, $x_1 = x$, $x_2 = -x$, $\lambda = \|x + z\|$, $x_3 \in C(p_4, \lambda)$ y $q \in C(O, \tfrac{\lambda}{2})$, entonces existen puntos $p_1$ y $p_2$ tales que $\{p_1, p_2, p_3, p_4\}$ y $\{x_1, x_2, x_3, x_4\}$ son sistemas $C$–ortocéntricos, donde $p_1 = S_q(x_1)$, $p_2 = S_q(x_2)$ y $x_4 = S_q(p_4)$. Además, haciendo $L_1 = \langle S_{x_1}(p_3), S_{x_2}(p_3) \rangle$, se cumple uno de los siguientes enunciados:

1. Si $p_1 = \dfrac{\lambda z}{\|z\|} - x$ y $p_2 = \dfrac{\lambda z}{\|z\|} + x$, $p_3$ y la línea $\langle p_1, p_2 \rangle$ están separados por $L_1$, o $L_1 = \langle p_1, p_2 \rangle$.
2. $\langle p_1, p_2 \rangle$ es la línea de soporte común de las circunferencias $C(x_2, \lambda)$ y $C(x_1, \lambda)$ si y solo si $q \perp_B x$.
3. Existen puntos $p_1$ y $p_2$ tal que $p_3$ y la línea $\langle p_1, p_2 \rangle$ están separados por $L_1$, o $L_1 = \langle p_1, p_2 \rangle$. Además, $p_4 \in A_G([p_3, p_1\rangle, [p_3, p_2\rangle)$.
4. Existen puntos $p_1$ y $p_2$ tal que $p_4 \in A_G^*([p_3, p_1\rangle, [p_3, p_2\rangle)$.

*Demostración:* Por el Lema 2.2 se tiene que existen puntos $p_1$ y $p_2$ tales que $\{p_1, p_2, p_3, p_4\}$ y $\{x_1, x_2, x_3, x_4\}$ son sistemas $C$-ortocéntricos. Probemos que existen puntos $p_1$ y $p_2$ que cumplen uno de los otros enunciados mencionados.

Tomemos los puntos $S_{x_i}(p_3) = d_i$ y $p_i = S_q(x_i)$, para $i = 1, 2$. Sean las rectas $L_2 = \langle p_1, p_2 \rangle$, $L_1 = \langle d_2, d_1 \rangle$ y $L_0 = \langle x_2, x_1 \rangle$, con vectores directores $p_2 - p_1$, $d_1 - d_2$ y $x_1 - x_2$, respectivamente. Como $d_1 - d_2 = 4x$, $p_2 - p_1 = 2x$ y $x_1 - x_2 = 2x$, se tiene que $L_0$, $L_1$ y $L_2$ son rectas paralelas, siendo $L_0$ exactamente el eje de las abscisas (ver *Figura 1*). También se cumple que $p_4 = \dfrac{d_1 + d_2}{2}$, $2q = \dfrac{p_1 + p_2}{2}$ y $O = \dfrac{x_1 + x_2}{2}$, de manera que $O \in L_0$, $p_4 \in L_1$ y $2q \in L_2$. Tomando la base $\{x, z\}$ de $M$, donde $x, z \neq O$ y tomando $q = \dfrac{\lambda z}{\|z\|}$, se calculan las coordenadas de $p_3$, $p_4$, $x_1$, $x_2$, $q$, $p_1$, $p_2$, $d_1$ y $d_2$, teniendo que $L_2 = \left(2t - 1, \dfrac{\lambda}{\|z\|}\right)$ y $L_1 = (4s - 1, 1)$, de manera que los puntos $p_4$ y $2q$ están montados en el eje de las ordenadas. Por tanto, las rectas $L_1$ y $L_2$ pasan por los puntos $(0, 1)$ y $\left(0, \dfrac{\lambda}{\|z\|}\right)$, respectivamente. Por otro lado,

$$\dfrac{2}{3}\|z\| = \dfrac{\|2z\|}{3} \leq \dfrac{\|z + x\|}{3} + \dfrac{\|z - x\|}{3}$$

pero como $\|z - x\| = \|z + x\|$, entonces $\|z\| \leq \lambda$. Así, en el caso donde $\dfrac{\lambda}{\|z\|} > 1$ se tiene que $L_1$ separa a la rectas $L_2$ del punto $p_3$ y si $\dfrac{\lambda}{\|z\|} = 1$, entonces $L_1 = L_2$.

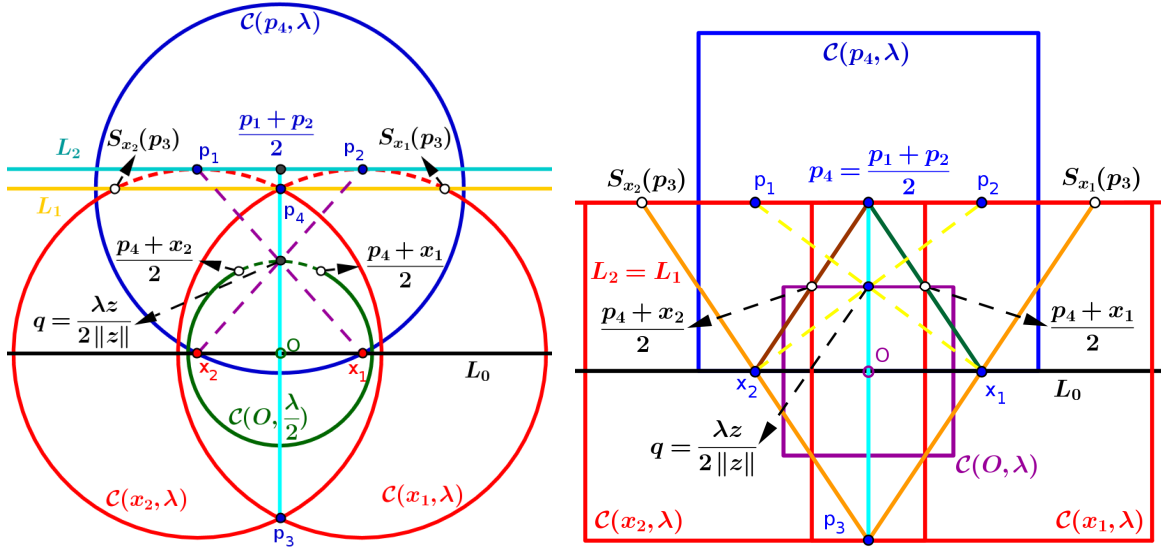

**Figura 1:** Demostración ítem 1 del Lema 3.1

2. Sea $q \in C\left(O, \dfrac{\lambda}{2}\right)$ tal que $q \perp_B x$, entonces la recta $Q$ que pasa por $q$, paralela a la línea $\langle x_1, x_2 \rangle$ (ver *Figura 2*), es la línea soporte de la circunferencia $C\left(O, \dfrac{\lambda}{2}\right)$. Aplicando las homotecias $H_{x_1,2}$ y $H_{x_2,2}$ a la recta $Q$ se obtienen las rectas $H_{x_1,2}(Q)$ y $H_{x_2,2}(Q)$, que pasan por $p_1 = 2q - x_1$ y $p_2 = 2q - x_2$, y son las líneas soporte de las circunferencias $C(x_2, \lambda)$ y $C(x_1, \lambda)$, respectivamente. Por otro lado, $\langle p_1, p_2 \rangle$ es paralela a $\langle x_1, x_2 \rangle$ por el ítem 2 del Teorema 2.1. Como por un punto pasa una sola recta paralela a otra recta dada, se tiene que $H_{x_1,2}(Q) = \langle p_1, p_2 \rangle = H_{x_2,2}(Q)$ y por tanto, $\langle p_1, p_2 \rangle$ es la línea soporte común de las circunferencias $C(x_1, \lambda)$ y $C(x_2, \lambda)$.

Recíprocamente, como $p_i = S_q(x_i)$, para $i = 1, 2$ por el ítem 2 del Teorema 2.1 se tiene que $x_1 - p_2 = x_2 - p_1$. Como $\langle p_1, p_2 \rangle$ es la línea soporte común de las circunferencias $C(x_1, \lambda)$ y $C(x_2, \lambda)$, entonces $x_1 - p_2 \perp_B x$ y $x_2 - p_1 \perp_B x$, de manera que $q \perp_B x$, por la homogeneidad de la ortogonalidad Birkhoff (ver *Figura 2*).

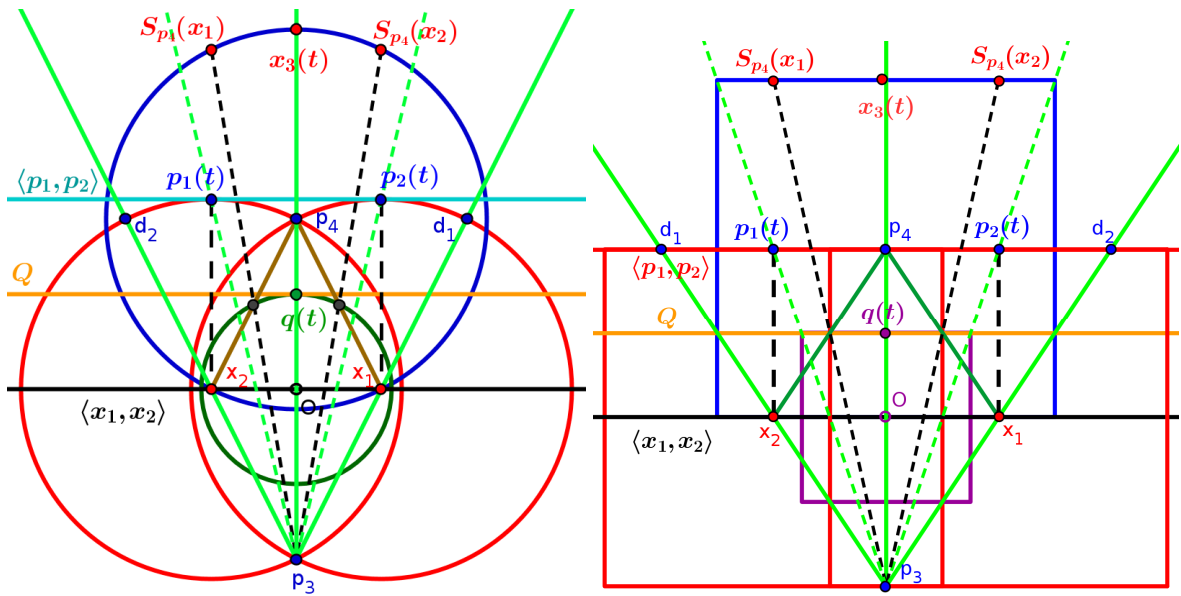

**Figura 2:** Demostración ítem 2 del Lema 3.1

3. Para cualquier $t \in [0,1]$, definamos las funciones continuas

$$x_3(t): [0,1] \to Arc^+_{C(p_4,\lambda)}\left(S_{p_4}(x_1), S_{p_4}(x_2)\right)$$

$q(t) = \dfrac{p_3 + x_3(t)}{2}$ y $p_i(t) = 2q(t) - x_i$ para $i = 1,2$. Cuando $t$ se mueve continuamente de $0$ a $1$, el rayo $A_G\left(\langle p_3, p_1(t)\rangle, \langle p_3, p_2(t)\rangle\right)$ se transforma continuamente de $A_G\left(\langle p_3, d_1\rangle, \langle p_3, p_4\rangle\right)$ a $A_G\left(\langle p_3, p_4\rangle, \langle p_3, d_2\rangle\right)$ (ver *Figura 1*). Así, existe $t_0 \in (0,1)$ tal que

$$A_G\left(\langle p_3, p_1(t_0)\rangle, \langle p_3, p_2(t_0)\rangle\right) = \langle p_3, p_4\rangle.$$

Sean $p_1 = p_1(t_0)$ y $p_2 = p_2(t_0)$, entonces $p_1$ y $p_2$ son puntos que tienen la propiedad deseada y, por el ítem 4 del Lema 2.1, $p_3$ y la línea $\langle p_1, p_2\rangle$ están separados por la línea $L_1$, que pasa por $p_4$ y es paralela a $\langle p_1, p_2\rangle$, o $L_1 = \langle p_1, p_2\rangle$ (ver *Figura 1*).

4. Definamos las funciones continuas $x_3(t)$, $q(t)$, $p_1(t)$ y $p_2(t)$ de igual forma que en el ítem 3 del presente lema. Tomemos la base del plano dada por $\Theta = \left\{\dfrac{x}{\|x\|}, \dfrac{z}{\|z\|}\right\}$ y sea $(a_3(t), b_3(t))$ el vector coordenado de $x_3(t)$ en dicha base para $t \in [0,1]$. Tomando en cuenta las coordenadas de los puntos $p_3$, $p_4$, $x_1$, $x_2$, , , en la base $\Theta$, se tiene que las correspondientes a $q(t)$, $p_1(t)$ y $p_2(t)$ son: $\left(\dfrac{a_3(t)}{2}, \dfrac{b_3(t)-\|z\|}{2}\right)$, $(a_3(t)-\|x\|, b_3(t)-\|z\|)$ y $(a_3(t)+\|x\|, b_3(t)-\|z\|)$, respectivamente. Como las rectas $\langle p_3, p_1(t) \rangle$ y $\langle p_3, p_2(t) \rangle$ están dadas por $kp_3 + (1-k)p_1(t)$ y $kp_3 + (1-k)p_2(t)$ con $k \in \mathbf{R}$, respectivamente, entonces estas se pueden expresar de la siguiente manera:

$$((1-k)(a_3(t)-\|x\|), -kb_3(t)+b_3(t)-\|z\|)$$

y

$$((1-k)(a_3(t)+\|x\|), -kb_3(t)+b_3(t)-\|z\|),$$

respectivamente. Deduciendo deducir que las ecuaciones de las rectas $\langle p_3, p_1(t) \rangle$ y $\langle p_3, p_2(t) \rangle$ en la base $\Theta$, están dadas por:

$$-b_3(t)x - (\|x\| - a_3(t))y + \|z\|(a_3(t) - \|x\|)$$

y

$$-b_3(t)x - (\|x\| + a_3(t))y + \|z\|(a_3(t) + \|x\|)$$

respectivamente. Por tanto, los elementos

$$f_{p_1(t)} = (-b_3(t), a_3(t) - \|x\|) \quad \text{y} \quad f_{p_2(t)} = (-b_3(t), -a_3(t) - \|x\|)$$

son funcionales asociados a las rectas $\langle p_3, p_1(t) \rangle$ y $\langle p_3, p_2(t) \rangle$, respectivamente. Como $A_G^*([p_3, p_1(t)], [p_3, p_2(t)])$ es un rayo que consta de los elementos $z = (w, e)$, tales que:

$$\dfrac{|f_{p_1(t)}(z)|}{\|f_{p_1(t)}\|} = \dfrac{|f_{p_2(t)}(z)|}{\|f_{p_2(t)}\|},$$

entonces los puntos $z$ anulan a la función continua definida por:

$$\phi_z(t) = \|f_{p_2(t)}\| \left| -wb_3(t) + e(a_3(t) - \|x\|) \right| - \|f_{p_1(t)}\| \left| -wb_3(t) - e(a_3(t) + \|x\|) \right|,$$

es decir,

$$A_G^*([p_3, p_1(t)], [p_3, p_2(t)]) = \{z \in M : \phi_z(t) = 0\}.$$

Fijando $p_4$, tenemos que

$$\phi_{p_4}(t) = \|f_{p_2(t)}\| \|z\| \left| \|x\| - a_3(t) \right| - \|f_{p_1(t)}\| \|z\| \left| \|x\| + a_3(t) \right|.$$

Como

$$\underset{t\to 0}{\text{Lím}}\, x_3(t) = S_{p_4}(x_1) = 2z - x \qquad y \qquad \underset{t\to 1}{\text{Lím}}\, x_3(t) = S_{p_4}(x_2) = 2z + x,$$

se tiene que

$$\underset{t\to 0}{\text{Lím}}\, (a_3(t), b_3(t)) = (-\|x\|, 2\|z\|) \qquad y \qquad \underset{t\to 1}{\text{Lím}}\, (a_3(t), b_3(t)) = (\|x\|, 2\|z\|)$$

Por tanto,

$$\underset{t\to 0}{\text{Lím}}\, \phi_{p_4}(t) = 2\|f_{p_2(t)}\|\|z\|\|x\| > 0 \qquad y \qquad \underset{t\to 1}{\text{Lím}}\, \phi_{p_4}(t) = -2\|f_{p_1(t)}\|\|z\|\|x\| < 0$$

de manera que existe un $t_0$ tal que $\phi_{p_4}(t_0) = 0$, es decir, existe $t_0$ tal que $p_4 \in A_G^*([p_3, p_1(t_0)\rangle, [p_3, p_2(t_0)\rangle)$. Así, los valores $p_1(t_0)$ y $p_2(t_0)$ son los puntos que cumplen la propiedad deseada. $\square$

Teniendo las propiedades de los sistemas $C$ – ortocéntricos, que se mostraron en el lema precedente, se probarán los resultados referentes a caracterización de euclidianidad en planos de Minkowski. El siguiente teorema dice que un plano de Minkowski es euclidiano si y solo si todo sistema $C$ – ortocéntrico $\{p_1, p_2, p_3, p_4\}$ con $p_3 \neq p_4$, cumple la siguiente implicación:

$$\begin{array}{c} \langle p_1, p_2 \rangle \text{ es la línea de soporte común de las} \\ \text{circunferencias que contienen los puntos} \\ \{p_1, p_3, p_4\} \text{ y } \{p_2, p_3, p_4\} \end{array} \Rightarrow \|p_3 - p_1\| = \|p_3 - p_2\|.$$

**Teorema 3.1:** Un plano de Minkowski $M$ es euclidiano si y solo si para cualquier sistema $C$ – ortocéntrico $\{p_1, p_2, p_3, p_4\}$ con $p_3 \neq p_4$, la igualdad $\|p_3 - p_1\| = \|p_3 - p_2\|$ se cumple, siempre y cuando $\langle p_1, p_2 \rangle$ es la línea de soporte común de las circunferencias que contienen los puntos $\{p_1, p_3, p_4\}$ y $\{p_2, p_3, p_4\}$, respectivamente.

*Demostración:* Supóngase $M$ euclídeo. Sea $\{p_1, p_2, p_3, p_4\}$ un sistema $C$ – ortocéntrico y $C(x_1, \lambda)$ la circunferencia circunscrita del triangulo $\Delta p_2 p_3 p_4$. Definamos $q = \dfrac{x_1 + p_1}{2}$ y $x_i = S_q(p_i)$ para $i = 2,3,4$, entonces $\langle p_3, p_4 \rangle$ es el eje radical de $C(x_1, \lambda)$ y $C(x_2, \lambda)$. Si $\langle p_1, p_2 \rangle$ es la línea soporte común de $C(x_1, \lambda)$ y $C(x_2, \lambda)$, entonces $\langle p_3, p_4 \rangle$ sería la bisectriz perpendicular de $\langle p_1, p_2 \rangle$ y por tanto $\|p_3 - p_1\| = \|p_3 - p_2\|$.

Recíprocamente, por el Lema 2.1 basta ver que para cualesquiera $x, z \in C$ con $x \perp_B z$ se cumple que $x \perp_I z$. Sea $x \in C$ y, sin perdida de generalidad, tomemos el arco $Arc_C^-(-x, x)$ y definamos la parametrización continua $\alpha : [0,1] \to Arc_C^-(-x, x)$ tal que $\alpha(0) = -x$ y $\alpha(1) = x$. Luego, por las propiedades de la ortogonalidad Birkhoff, existe $z \in Arc_C^-(-x, x)$ tal que $z \perp_B x$ y por tanto, hay $k_1 \in [0,1]$ tal que $\alpha(k_1) = z$.

Sea $J = \{h \in [0,1] : x \perp_I \alpha(h)\}$ y tomemos $h_1 \in J$ tal que $|h_1 - k_1| = \min_{h \in J}\{|h - k_1|\}$. Ahora, para cualquier $t > 0$, $tx$ y $-tx$ son los puntos de intersección de la circunferencia $tC = C(O, t)$ y la recta $\langle -x, x \rangle$. Sea la parametrización continua $\alpha_t : [0,1] \to Arc_{tC}^-(-tx, tx)$, con $\alpha_t(0) = -tx$ y $\alpha_t(1) = tx$. Sea $z_t \in Arc_{tC}^-(-tx, tx)$ tal que $x \perp_I z_t$, por tanto existe $k_t \in [0,1]$ tal que $\alpha_t(k_t) = z_t$. Sea $J_t = \{h \in [0,1] : x \perp_I \alpha_t(h)\}$ y definamos la función continua $y(t) = \alpha_t(h_t)$ con $h_t \in J_t$ tal que $|h_t - k_t| = \min_{h \in J}\{|h - k_t|\}$.

Hagamos $p_4(t) = y(t)$, $p_3(t) = -y(t)$, $x_1 = x$ y $x_2 = -x$. Tomemos $q(t) = \dfrac{x + y(t)}{2}$ y definamos los puntos $p_1(t) = S_{q(t)}(x_1)$ y $p_2(t) = S_{q(t)}(x_2)$. Por el ítem 2 del Lema 3.1, $\langle p_1(t), p_2(t) \rangle$ es la línea de soporte común de las circunferencias que contienen los puntos $\{p_1(t), p_3(t), p_4(t)\}$ y $\{p_2(t), p_3(t), p_4(t)\}$, respectivamente. Por el ítem 2 del Teorema 2.1 se tiene que $p_1(t) - p_2(t) = x_2 - x_1$. Así, $p_1(t) - x_2 \perp_B x$ y $p_2(t) - x_1 \perp_B x$, de manera que:

$$p_1(t) - x_2 = p_2(t) - x_1 = 2q(t) = \|x + y(y)\| z.$$

Sea $z(t) = \|x + y(t)\| z$, entonces

$$\|p_3(t) - p_1(t)\| = \|p_3(t) - x_2 + x_2 - p_1(t)\| = \|p_3(t) - x_2 - z(t)\|$$

y

$$\|p_3(t) - p_1(t)\| = \|p_3(t) - x_1 + x_1 - p_1(t)\| = \|p_3(t) - x_1 - z(t)\|.$$

Por hipótesis, $\|p_3 - p_1(t)\| = \|p_3 - p_2(t)\|$ y por tanto,

$$\|p_3(t) - x_2 - z(t)\| = \|p_3(t) - x_1 - z(t)\|,$$

es decir, $\|(y(t) - z(t)) + x\| = \|(y(t) - z(t)) - x\|$. Es claro que $\lim_{t \to 0} y(t) = O$ y por tanto, $\lim_{t \to 0} \|x + y(t)\| = \|x\| = 1$. De esto se sigue que

$$\text{Lím}_{t \to 0} (y(t) - z(t)) = \text{Lím}_{t \to 0} (y(t) - \|x + y(t)\| z) = -z,$$

y entonces

$$\|x - z\| = \left\|\text{Lím}_{t \to 0} (x + y(t) + z(t))\right\| = \left\|\text{Lím}_{t \to 0} (y(t) - z(t) - x)\right\| = \|x + z\|,$$

de manera que $z \perp_I x$. □

El siguiente resultado dice que un plano de Minkowski es euclidiano si y solo si todo sistema $C$–ortocéntrico $\{p_1, p_2, p_3, p_4\}$ con $p_3 \neq p_4$, cumple la siguiente implicación:

$$\|p_3 - p_1\| = \|p_3 - p_2\| \quad \Rightarrow \quad \begin{array}{l} \langle p_1, p_2 \rangle \text{ es la línea de soporte común de las} \\ \text{circunferencias que contienen los puntos} \\ \{p_1, p_3, p_4\} \text{ y } \{p_2, p_3, p_4\} \end{array}$$

**Teorema 3.2:** Un plano de Minkowski $M$ es euclidiano si y solo si para cualquier sistema $C$–ortocéntrico $\{p_1, p_2, p_3, p_4\}$ con $p_3 \neq p_4$, $\langle p_1, p_2 \rangle$ es la línea de soporte común de las circunferencias que contienen los puntos $\{p_1, p_3, p_4\}$ y $\{p_2, p_3, p_4\}$, respectivamente, siempre y cuando $\|p_3 - p_1\| = \|p_3 - p_2\|$.

*Demostración:* La primera implicación se deduce empleando un argumento similar al usado en el Teorema 3.1.

Recíprocamente, para cualquier $x, z \in C$ con $z \perp_B x$ y cualquier $t > 0$, definamos $y(t)$, $x_1$, $x_2$, $p_3(t)$ y $p_4(t)$ como en la prueba del Teorema 3.1. Entonces, por el ítem 2 del Lema 2.2, existen puntos $p_1(t)$ y $p_2(t)$ tales que el conjunto $\{p_1(t), p_2(t), p_3(t), p_4(t)\}$ es un sistema $C$–ortocéntrico con $\|p_3(t) - p_1(t)\| = \|p_3(t) - p_2(t)\|$. Por el ítem 1 del Lema 3.1, $p_3(t)$ y la línea $\langle p_1(t), p_2(t) \rangle$ están separados por la recta que pasa por $p_4(t)$, paralela a $\langle p_1(t), p_2(t) \rangle$, o $\langle p_1(t), p_2(t) \rangle = L_1$. Por hipótesis, $\langle p_1(t), p_2(t) \rangle$ es la línea de soporte común de las circunferencias $C(x_1, \|x + y(t)\|)$ y $C(x_2, \|x + y(t)\|)$. Luego, como en la prueba del Teorema 3.1, se puede demostrar que $z \perp_I x$. □

El siguiente enunciado dice que un plano de Minkowski es euclidiano si y solo si todo sistema $C$–ortocéntrico $\{p_1, p_2, p_3, p_4\}$ con $p_3 \neq p_4$, cumple la siguiente implicación:

$$p_4 \in \left\langle p_3, \frac{p_1 + p_2}{2} \right\rangle \quad \Rightarrow \quad \begin{array}{l} \langle p_1, p_2 \rangle \text{ es la línea de soporte común de las} \\ \text{circunferencias que contienen los puntos} \\ \{p_1, p_3, p_4\} \text{ y } \{p_2, p_3, p_4\} \end{array}$$

**Teorema 3.3:** Un plano de Minkowski $M$ es euclidiano si y solo si para cualquier sistema $C$–ortocéntrico $\{p_1, p_2, p_3, p_4\}$ con $p_3 \neq p_4$, $\langle p_1, p_2 \rangle$ es la línea de soporte común de las circunferencias que contienen los puntos $\{p_1, p_3, p_4\}$ y $\{p_2, p_3, p_4\}$, respectivamente, siempre y cuando $p_4 \in \left\langle p_3, \frac{p_1 + p_2}{2} \right\rangle$.

*Demostración:* Sea $\{p_1, p_2, p_3, p_4\}$ un sistema $C$–ortocéntrico y $C(x_1, \lambda)$ la circunferencia circunscrita del triangulo $\Delta p_2 p_3 p_4$. Definamos $q = \frac{x_1 + p_1}{2}$ y $x_i = S_q(p_i)$ para $i = 2, 3, 4$. Por el ítem 3 del Teorema 2.1 se tiene que el cuadrilátero formado por los puntos $x_1$, $p_3$, $x_2$, $p_4$ es un rombo y, como el plano es euclídeo, entonces $\langle p_3, p_4 \rangle$ y $\langle x_1, x_2 \rangle$ son perpendiculares y se cortan en el punto $\frac{x_1 + x_2}{2}$. Como $p_4 \in \left\langle p_3, \frac{p_1 + p_2}{2} \right\rangle$ se tiene que el segmento $\left[ \frac{p_1 + p_2}{2}, \frac{x_1 + x_2}{2} \right]$ es perpendicular a la recta $\langle x_1, x_2 \rangle$ en el mismo punto $\frac{x_1 + x_2}{2}$.

Por el ítem 2 del Teorema 2.1, el segmento $\left[ \frac{p_1 + p_2}{2}, \frac{x_1 + x_2}{2} \right]$ es paralelo a los segmentos $[p_1, x_2]$ y $[p_2, x_1]$. Además, estos últimos son perpendiculares a $\langle x_1, x_2 \rangle$ en los puntos $x_2$ y $x_1$, respectivamente. Así, la recta $L_1$ paralela a $\langle x_1, x_2 \rangle$, que pasa por $p_1$, es de soporte de $C(x_2, \lambda)$. De igual forma, la recta $L_2$ paralela a $\langle x_1, x_2 \rangle$, que pasa por $p_2$, es de soporte de $C(x_1, \lambda)$, pero como $\langle p_1, p_2 \rangle$ es paralela a $\langle x_1, x_2 \rangle$ entonces $L_1 = L_2 = \langle p_1, p_2 \rangle$.

Recíprocamente, sean $x, z \in M$ tal que $x \perp_I z$. Hagamos $p_4 = z$, $p_3 = -z$, $x_1 = x$, $x_2 = -x$. Por el ítem 3 del Lema 2.2, se tiene que existe un sistema $C$–ortocéntrico

$\{p_1, p_2, p_3, p_4\}$ con $p_3 \neq p_4$, tal que $p_4 \in \left\langle p_3, \dfrac{p_1 + p_2}{2} \right\rangle$ y por tanto, se tiene que $p_1 = \pm \dfrac{\lambda z}{\|z\|} - x$ y $p_2 = \pm \dfrac{\lambda z}{\|z\|} + x$.

Por hipótesis, $\langle p_1, p_2 \rangle$ es la línea de soporte común de las circunferencias que contienen a los puntos $\{p_1, p_3, p_4\}$ y $\{p_2, p_3, p_4\}$, respectivamente. Como las rectas $\langle p_1, p_2 \rangle$ y $\langle x_1, x_2 \rangle$ son paralelas, entonces $p_1 - x_2 \perp_B x$ y $p_2 - x_1 \perp_B x$. Por lo tanto $\pm \dfrac{\lambda z}{\|z\|} \perp_B x$, y por la homogeneidad de la ortogonalidad Birkhoff, se tiene que $z \perp_B x$. Luego, por el Lema 2.1 $M$ es euclídeo. $\square$

El siguiente teorema dice que un plano de Minkowski es euclidiano si y solo si todo sistema $C$–ortocéntrico $\{p_1, p_2, p_3, p_4\}$ con $p_3 \neq p_4$, cumple la siguiente implicación:

$\langle p_1, p_2 \rangle$ es la línea de soporte común de las
circunferencias que contienen los puntos $\quad \Rightarrow \quad p_4 \in \left\langle p_3, \dfrac{p_1 + p_2}{2} \right\rangle$
$\{p_1, p_3, p_4\}$ y $\{p_2, p_3, p_4\}$

**Teorema 3.4:** Un plano de Minkowski $M$ es euclidiano si y solo si para cualquier sistema $C$–ortocéntrico $\{p_1, p_2, p_3, p_4\}$ con $p_3 \neq p_4$, se cumple que $p_4 \in \left\langle p_3, \dfrac{p_1 + p_2}{2} \right\rangle$, siempre y cuando $\langle p_1, p_2 \rangle$ es la línea de soporte común de las circunferencias que contienen los puntos $\{p_1, p_3, p_4\}$ y $\{p_2, p_3, p_4\}$, respectivamente.

*Demostración:* Sea un sistema $C$–ortocéntrico $\{p_1, p_2, p_3, p_4\}$. Usando el mismo argumento del Teorema 3.3 se tiene que $\left\langle p_4, \dfrac{x_1 + x_2}{2} \right\rangle$ es perpendicular a $\langle x_1, x_2 \rangle$ en el punto $\dfrac{x_1 + x_2}{2}$. Como $\langle p_1, p_2 \rangle$ es la línea de soporte común de las circunferencias $C(x_1, \lambda)$ y $C(x_2, \lambda)$, se tiene que los segmentos $[x_1, p_2]$ y $[x_2, p_1]$ son perpendiculares a los puntos $\{x_1, x_2, p_1, p_2\}$ es un paralelogramo y por tanto el segmento $\left[ \dfrac{p_1 + p_2}{2}, \dfrac{x_1 + x_2}{2} \right]$, paralelo a $[x_1, p_2]$, es perpendicular a $\langle x_1, x_2 \rangle$ en el punto $\dfrac{x_1 + x_2}{2}$, y por la unicidad de la perpendicularidad (salvo el opuesto) se tiene:

$$\left\langle \frac{p_1 + p_2}{2}, \frac{x_1 + x_2}{2} \right\rangle \subset \langle p_3, p_4 \rangle$$

de manera que $p_4 \in \left\langle p_3, \frac{p_1 + p_2}{2} \right\rangle$.

Recíprocamente, sean $x, z \in M$ tal que $x \perp_I z$. Hagamos $p_4 = z$, $p_3 = -z$, $x_1 = x$, $x_2 = -x$. Por el ítem 2 del Lema 2.2, se tiene que existe un sistema $C$ – ortocéntrico $\{p_1, p_2, p_3, p_4\}$ con $p_3 \neq p_4$, tal que $\langle p_1, p_2 \rangle$ es la línea de soporte común de las circunferencias que contienen los puntos $\{p_1, p_3, p_4\}$ y $\{p_2, p_3, p_4\}$, respectivamente. Como las rectas $\langle p_1, p_2 \rangle$ y $\langle x_1, x_2 \rangle$ son paralelas, entonces $p_1 - x_2 \perp_B x$ y $p_2 - x_1 \perp_B x$. Como $p_4 \in \left\langle p_3, \frac{p_1 + p_2}{2} \right\rangle$, entonces por el ítem 2 del Lema 2.2 se tiene que $p_1 = \pm \frac{\lambda z}{\|z\|} - x$ y $p_2 = \pm \frac{\lambda z}{\|z\|} + x$. Así, $\pm \frac{\lambda z}{\|z\|} \perp_B x$ de manera $z \perp_B x$ y por tanto, $M$ es euclídeo. $\square$

El siguiente resultado dice que un plano de Minkowski es euclidiano si y solo si todo sistema $C$ – ortocéntrico $\{p_1, p_2, p_3, p_4\}$ con $p_3 \neq p_4$, cumple la siguiente implicación:

$$p_4 \in \left\langle p_3, \frac{p_1 + p_2}{2} \right\rangle \implies p_4 \in A_B([p_3, p_1), [p_3, p_2))$$

**Teorema 3.5:** Un plano de Minkowski $M$ es euclidiano si y solo si para cualquier sistema $C$ – ortocéntrico $\{p_1, p_2, p_3, p_4\}$ con $p_3 \neq p_4$, se cumple que $p_4 \in A_B([p_3, p_1), [p_3, p_2))$, siempre y cuando $p_4 \in \left\langle p_3, \frac{p_1 + p_2}{2} \right\rangle$.

*Demostración:* Sea $M$ un plano euclidiano y $\{p_1, p_2, p_3, p_4\}$ un sistema $C$ – ortocéntrico $\{p_1, p_2, p_3, p_4\}$ con $p_3 \neq p_4$. Como en un plano euclídeo la bisectriz de Busemann coincide con la bisectriz euclidiana, entonces $A_B([p_3, p_1), [p_3, p_2)) = \left[ p_3, \frac{p_1 + p_2}{2} \right\rangle$ y como $p_4 \in \left\langle p_3, \frac{p_1 + p_2}{2} \right\rangle$, se tiene lo deseado.

Recíprocamente, sean $x, z \in M$ tal que $x \perp_I z$. Hagamos $p_4 = z$, $p_3 = -z$, $x_1 = x$, $x_2 = -x$. Por el ítem 3 del Lema 2.2, existe un sistema $C$-ortocéntrico $\{p_1, p_2, p_3, p_4\}$ con $p_3 \neq p_4$, tal que $p_4 \in \left\langle p_3, \dfrac{p_1 + p_2}{2} \right\rangle$, entonces $p_4 \in A_B([p_3, p_1\rangle, [p_3, p_2\rangle)$, por hipótesis. Luego, por el ítem 3 del Lema 2.2, $p_1 = \pm\dfrac{\lambda z}{\|z\|} - x$ y $p_2 = \pm\dfrac{\lambda z}{\|z\|} + x$. Tomando los valores $p_1 = \dfrac{\lambda z}{\|z\|} - x$ y $p_2 = \dfrac{\lambda z}{\|z\|} + x$, entonces $\dfrac{p_1 + p_2}{2} = \dfrac{\lambda z}{\|z\|}$.

Por otra parte, existen $\alpha \in \mathbf{R}$ y $\beta \in [1, +\infty)$ tales que

$$p_4 = \alpha\, p_3 + (1-\alpha)\dfrac{p_1 + p_2}{2} \quad \text{y} \quad p_4 = \beta\, p_3 + (1-\beta)\left(p_3 + \dfrac{p_3 - p_1}{\|p_3 - p_1\|} + \dfrac{p_3 - p_2}{\|p_3 - p_2\|}\right)$$

y sustituyendo los valores de los puntos $p_1$, $p_2$ y $p_3$ se tiene que

$$p_4 = \left(\dfrac{\lambda - \alpha(\|z\| + \lambda)}{\|z\|}\right) z \quad \text{y} \quad p_4 = rz + (1-\beta)\left(\dfrac{1}{\left\|x - \left(\dfrac{\|z\| + \lambda}{\|z\|}\right)z\right\|} + \dfrac{1}{\left\|x + \left(\dfrac{\|z\| + \lambda}{\|z\|}\right)z\right\|}\right) x$$

respectivamente, donde $r = (\beta - 1)\left(1 + \dfrac{\|z\| + \lambda}{\|z\|\left\|x - \left(\dfrac{\|z\| + \lambda}{\|z\|}\right)z\right\|} + \dfrac{\|z\| + \lambda}{\|z\|\left\|x + \left(\dfrac{\|z\| + \lambda}{\|z\|}\right)z\right\|}\right) - 1$. Por tanto,

$$\left\|x - \left(\dfrac{\|z\| + \lambda}{\|z\|}\right)z\right\| = \left\|x + \left(\dfrac{\|z\| + \lambda}{\|z\|}\right)z\right\|$$

Haciendo $t = \dfrac{\|z\| + \lambda}{\|z\|}$ se tiene que $x \perp_I tz$ y así $M$ es euclidiano. □

El siguiente enunciado dice que un plano de Minkowski es euclidiano si y solo si todo sistema $C$-ortocéntrico $\{p_1, p_2, p_3, p_4\}$ con $p_3 \neq p_4$, cumple la siguiente implicación:

$$p_4 \in A_B(\langle p_3, p_1 \rangle, \langle p_3, p_2 \rangle) \implies p_4 \in \left\langle p_3, \frac{p_1 + p_2}{2} \right\rangle$$

**Teorema 3.6:** Un plano de Minkowski $M$ es euclidiano si y solo si para cualquier sistema $C$ – ortocéntrico $\{p_1, p_2, p_3, p_4\}$ con $p_3 \neq p_4$, se cumple que $p_4 \in \left\langle p_3, \frac{p_1 + p_2}{2} \right\rangle$, siempre y cuando $p_4 \in A_B(\langle p_3, p_1 \rangle, \langle p_3, p_2 \rangle)$.

**Demostración:** Si $M$ es euclídeo el enunciado del teorema es claramente cierto. Recíprocamente, sean $x, z \in M$ tal que $x \perp_I z$. Hagamos $p_4 = z$, $p_3 = -z$, $x_1 = x$, $x_2 = -x$. Por el ítem 5 del Lema 2.2, existe un sistema $C$ – ortocéntrico $\{p_1, p_2, p_3, p_4\}$ con $p_3 \neq p_4$, tal que $p_4 \in A_B(\langle p_3, p_1 \rangle, \langle p_3, p_2 \rangle)$, entonces por hipótesis, $p_4 \in \left\langle p_3, \frac{p_1 + p_2}{2} \right\rangle$. Empleando el mismo argumento presentado en el Teorema 3.5, se tiene que $x \perp_I tz$ para $t = \frac{\|z\| + \lambda}{\|z\|}$ y por tanto, $M$ es euclidiano. □

El siguiente teorema dice que un plano de Minkowski es euclidiano si y solo si todo sistema $C$ – ortocéntrico $\{p_1, p_2, p_3, p_4\}$ con $p_3 \neq p_4$, cumple la siguiente implicación:

$$p_4 \in \left\langle p_3, \frac{p_1 + p_2}{2} \right\rangle \implies p_4 \in A_G^*(\langle p_3, p_1 \rangle, \langle p_3, p_2 \rangle)$$

**Teorema 3.7:** Un plano de Minkowski $M$ es euclidiano si y solo si para cualquier sistema $C$ – ortocéntrico $\{p_1, p_2, p_3, p_4\}$ con $p_3 \neq p_4$, se cumple que $p_4 \in A_G^*(\langle p_3, p_1 \rangle, \langle p_3, p_2 \rangle)$, siempre y cuando $p_4 \in \left\langle p_3, \frac{p_1 + p_2}{2} \right\rangle$.

*Demostración:* Si el plano $M$ es euclídeo, todo sistema $C$ – ortocéntrico $\{p_1, p_2, p_3, p_4\}$ con $p_3 \neq p_4$, satisface que $p_4$ está en la bisectriz euclidiana, siempre y cuando $p_4 \in \left\langle p_3, \frac{p_1 + p_2}{2} \right\rangle$ y por tanto, $d(p_4, L_{p_1}) = d(p_4, L_{p_2})$ donde $L_{p_1}$ y $L_{p_2}$ son las parametrizaciones de las rectas $\langle p_3, p_1 \rangle$ y $\langle p_3, p_2 \rangle$, con $p_3 - p_1$ y $p_3 - p_2$ sus vectores directores, respectivamente. Por tanto,

$$\frac{|L_{p_1}(p_4)|}{\|p_3-p_1\|^*} = \frac{|L_{p_2}(p_4)|}{\|p_3-p_2\|^*},$$

pero como en un plano euclídeo $\|\circ\|^* = \|\circ\|$, entonces

$$\frac{|L_{p_1}(p_4)|}{\|p_3-p_1\|} = \frac{|L_{p_2}(p_4)|}{\|p_3-p_2\|},$$

de manera que $p_4 \in A_G^*(\langle p_3,p_1\rangle,\langle p_3,p_2\rangle)$, teniendo así lo deseado.

Recíprocamente, sean $x,z \in M$ tales que $x \perp_I z$. Por el Lema 2.1 basta ver que existe un $t \in \mathbf{R}^+ - \{1\}$ tal que $z \perp_I tx$. Por el ítem 3 del Lema 2.2 existe un sistema $C$-ortocéntrico $\{p_1,p_2,p_3,p_4\}$ tal que $p_4 \in \left\langle p_3, \dfrac{p_1+p_2}{2}\right\rangle$, por tanto se pueden tomar $p_1 = \dfrac{\lambda z}{\|z\|} - x$ y $p_2 = \dfrac{\lambda z}{\|z\|} + x$.

Por otro lado, el conjunto $\Theta = \left\{\dfrac{x}{\|x\|}, \dfrac{z}{\|z\|}\right\}$ es una base de $M$, pues $x \perp_I z$. De manera que: $x_1 = (\|x\|,0)$, $x_2 = (-\|x\|,0)$, $p_3 = (0,-\|z\|)$, $p_4 = (0,\|z\|)$, $p_1 = (-\|x\|,\lambda)$ y $p_2 = (\|x\|,\lambda)$. Luego, aplicando la misma idea usada en el ítem 4 del Lema 3.1, se tiene que las ecuaciones de las rectas $\langle p_3,p_1\rangle$ y $\langle p_3,p_2\rangle$, en la base $\Theta$, están dadas por:

$$-(\|z\|+\lambda)x - \|x\|z - \|x\|\|z\| \qquad \text{y} \qquad -(\|z\|+\lambda)x + \|x\|z + \|x\|\|z\|$$

respectivamente. Así, los elementos

$$f_{p_1} = (\lambda + \|z\|, \|x\|) \qquad \text{y} \qquad f_{p_2} = (\lambda + \|z\|, -\|x\|)$$

son funcionales asociados a dichas rectas, respectivamente. Como por hipótesis $p_4 \in A_G^*(\langle p_3,p_1\rangle,\langle p_3,p_2\rangle)$, entonces $d(p_4,f_{p_1}) = d(p_4,f_{p_2})$. Por tanto,

$$\frac{|f_{p_1}(p_4)|}{\|f_{p_1}\|} = \frac{|f_{p_2}(p_4)|}{\|f_{p_2}\|}$$

de manera que

$$\frac{\|x\|\|z\|}{\left\|\frac{(\lambda+\|z\|)x}{\|x\|}+\frac{\|x\|z}{\|z\|}\right\|}=\frac{\|x\|\|z\|}{\left\|\frac{(\lambda+\|z\|)x}{\|x\|}-\frac{\|x\|z}{\|z\|}\right\|},$$

de donde

$$\left\|(\lambda+\|z\|)\|z\|x-\|x\|^2 z\right\|=\left\|(\lambda+\|z\|)\|z\|x+\|x\|^2 z\right\|,$$

teniendo que

$$\left\|\frac{(\lambda+\|z\|)\|z\|x}{\|x\|^2}-z\right\|=\left\|\frac{(\lambda+\|z\|)\|z\|x}{\|x\|^2}+z\right\|,$$

Tomando $t=\dfrac{(\lambda+\|z\|)\|z\|}{\|x\|^2}>0$ se obtiene que $z\perp_I tx$.

El siguiente resultado dice que un plano de Minkowski es euclidiano si y solo si todo sistema $C$–ortocéntrico $\{p_1,p_2,p_3,p_4\}$ con $p_3\neq p_4$, cumple la siguiente implicación:

$$p_4\in A_G^*([p_3,p_1\rangle,[p_3,p_2\rangle)\;\Rightarrow\;p_4\in\left\langle p_3,\frac{p_1+p_2}{2}\right\rangle$$

**Teorema 3.8:** Un plano de Minkowski $M$ es euclidiano si y solo si para cualquier sistema $C$–ortocéntrico $\{p_1,p_2,p_3,p_4\}$ con $p_3\neq p_4$, se cumple que $p_4\in\left\langle p_3,\dfrac{p_1+p_2}{2}\right\rangle$ siempre y cuando $p_4\in A_G^*([p_3,p_1\rangle,[p_3,p_2\rangle)$.

*Demostración:* Sea $\{p_1,p_2,p_3,p_4\}$ un sistema $C$-ortocéntrico en un plano euclídeo, con $p_3\neq p_4$, donde $p_4\in A_G^*([p_3,p_1\rangle,[p_3,p_2\rangle)$. Usando el razonamiento expuesto en el Teorema 3.7 se tiene que $p_4\in A_G([p_3,p_1\rangle,[p_3,p_2\rangle)$ y por tanto, $p_4$ está en la bisectriz euclidiana, de manera que $p_4\in\left\langle p_3,\dfrac{p_1+p_2}{2}\right\rangle$.

Recíprocamente, sean $x,z\in M$ tales que $x\perp_I z$. Por el Lema 2.1 basta ver que existe un $t\in\mathbf{R}^+-\{1\}$ tal que $z\perp_I tx$. Por el ítem 4 del Lema 3.1 existe un sistema

$C$–ortocéntrico $\{p_1,p_2,p_3,p_4\}$ tal que $p_4 \in A_G^*\bigl([p_3,p_1),[p_3,p_2)\bigr)$. Por hipótesis se tiene que $p_4 \in \left\langle p_3, \dfrac{p_1+p_2}{2} \right\rangle$ y, por el ítem 2 del Lema 2.1, se puede tomar $p_1 = \dfrac{\lambda z}{\|z\|} - x$ y $p_2 = \dfrac{\lambda z}{\|z\|} + x$. Usando la misma idea expuesta en el Teorema 3.7 se tiene que

$$\left\| \dfrac{(\lambda+\|z\|)\|z\|x}{\|x\|^2} - z \right\| = \left\| \dfrac{(\lambda+\|z\|)\|z\|x}{\|x\|^2} + z \right\|,$$

Tomando $t = \dfrac{(\lambda+\|z\|)\|z\|}{\|x\|^2} > 0$ se obtiene que $z \perp_I tx$.

El siguiente teorema relaciona la euclidianidad de un plano de Minkowski con la noción de ortogonalidad cordal y los sistemas $C$-ortocéntricos.

**Teorema 3.9:** Sean $M$ un plano de Minkowski y $[p_1,q_1]$, $[p_2,q_2]$ cuerdas en $C$, tales que $[p_1,q_1] \perp_C [p_2,q_2]$. Sea $h$ el $C$–ortocentro del $\Delta\, p_1q_1p_2$ asociado al origen. Si $\langle p_2, h \rangle \cap C = \{q_2\}$, entonces $M$ es euclidiano.

*Demostración:* Sean $[p_1,q_1]$ y $[p_2,q_2]$ dos cuerdas de $C$ tales que $[p_1,q_1] \perp_C [p_2,q_2]$, entonces lo segmentos $[p_1,q_1]$ y $[p_2,-q_2]$ son cuerdas paralelas de $C$. Como por hipótesis $h$ es el $C$–ortocentro del $\Delta\, p_1q_1p_2$ y $\langle p_2, h \rangle \cap C = \{q_2\}$, entonces $p_2$, $h$ y $q_2$ son colineales, de manera que $h = \alpha\, p_2 + (1-\alpha)(q_2)$ con $\alpha \in \mathbf{R}$.

Como el $\Delta\, p_1q_1p_2$ está inscrito en $C$, se tiene que $h = p_1 + q_1 + p_2$. Así,

$$\alpha\, p_2 + (1-\alpha)(q_2) = p_1 + q_1 + p_2$$

por tanto $(\alpha-1)p_2 + (\alpha-1)(-q_2) = p_1 + q_1$, de donde $(\alpha-1)\dfrac{(p_2-q_2)}{2} = \dfrac{(p_1+q_1)}{2}$. De esta forma, los puntos medios de los segmentos $[p_1,q_1]$ y $[p_2,-q_2]$ están alineados con $O$, y por el Lema 2.1 se tiene que $M$ es euclidiano.

### Referencias

[1] A. C. Thompson (1996). *Minkowski geometry*. Encyclopedia of Mathematics and Its Applications. **63**. Cambridge University Press. Cambridge. ISBN 0-521-40472-X.

[2] D. Amir (1986). *Characterizations of inner product spaces.* Birkhüser. Basel. ISBN 3-7643-1774-4.

[3] H. Busemann (1975). Planes with analogues to euclidean angular bisectors. *Math. Scand.* **36**, 5-11.

[4] G. Birkhoff (1935). Orthogonality in linear metric spaces. Duke Math. J. **1**(2), 169-172.

[5] H. Martini and K. J. Swanepoel (2006). Antinorms and radon curves. *Aequationes Math.* **72**, 110-138.

[6] H. Martini, K. J. Swanepoel and G. Weiβ (2001). The geometry of Minkowski spaces - A survey. Part I. *Expositiones Math.* **19**, 97-142.

[7] H. Martini and K. J. Swanepoel (2004). The geometry of Minkowski spaces - A survey. Part II. *Expositiones Math.* **22**, 93-144.

[8] H. Martini and M. Spirova (2010). A new type of orthogonality for normed planes. *Czechoslovak Mathematical Journal.* **60**(2), 339-349.

[9] H. Martini and S. Wu (2009). On orthocentric systems in strictly convex normed planes. *Extracta Mathematicae.* **24**(1), 31-45.

[10] J. Alonso, H. Martini and S. Wu (2012). On Birkhoff ortghogonality and isosceles orthogonality in normed linear spaces. *Aequationes Math.* **83**, 153-189.

[11] J. Alonso (1994). Uniqueness properties of isosceles orthogonality in normed linear spaces. *Ann. Sci. Math. Québec.* **18**(1), 25-38.

[12] R. A. Johson (2007). *Advanced euclidean geometry.* Dover Publications, Inc., Mineola, New York. ISBN-10: 0-486-46237-4.

[13] R. C. James (1945). Orthogonality in normed linear spaces. *Duke Math. J.* **12**, 291-302.

[14] S. Wu (2009). *Geometry of Minkowski planes and spaces - Selected topics.* Ph. D. Thesis. Chemnitz University of Technology.

[15] Rosas T. y Pacheco W. (2014). Orthocentric systems in Minkowski planes. *Beiträge zur Algebra und Geometrie (BZAG).* DOI 10.1007/s13366-014-0214-6. ISSN 0138-4821.

[16] Rosas T. (2014). *Sistemas Ortocéntricos en planos de Minkowski y euclidianidad.* Tesis Doctoral. Universidad Centroccidental Linsandro Alvarado, Barquisimeto, Venezuela.

[17] V. V. Glogovskij (1970). Bisectors on the Minkowski plane with norm $(x^p + y^p)^{\frac{1}{p}}$. *Visnik L'viv. Politehn. Inst.* **218**, 192-198.